\documentstyle[amssymb,amsmath,amsthm,amscd,graphicx,12pt]{article}

\title{Replacement of fixed sets for compact group actions:  The $2\rho$ theorem}
\author{Sylvain Cappell\thanks{Research was partially supported by a DARPA grant} \\
Courant Institute of Mathematical Sciences, New York University \\
\\
Shmuel Weinberger\thanks{Research was partially supported by NSF grants DMS 0504721 and 0805913} \\
University of Chicago \\
\\
Min Yan\thanks{Research was supported by Hong Kong RGC General Research Fund 605005 and 604408} \\ 
Hong Kong University of Science and Technology}
\DeclareMathOperator{\proj}{proj}

\newcommand{\sub}{\subset}

\newcommand{\pa}{\partial}

\newcommand{\bb}{\mathbb}
\newcommand{\mc}{\mathcal}
\newcommand{\s}{\hspace{2mm}}

\theoremstyle{definition}

\theoremstyle{remark}

\numberwithin{equation}{section}

\begin{document}
\maketitle

%\begin{abstract}
%\end{abstract}

\section{Introduction and Main Results}

The replacement problem for stratified spaces asks whether any manifold (simple) homotopy equivalent to the bottom stratum of a stratified space $X$ is the bottom stratum of a stratified space stratified (simple) homotopy equivalent to $X$.  
 
This is impossible in general, as one sees by considering $X = M \times [0,1]$ and working rel $M \times 0$. On the other hand, the classic theorem of Browder, Casson, Haefliger, Sullivan and Wall (see \cite{Wa1} Corollary 11.3.1) asserts that if $X = (W, M)$ is a pair consisting of a manifold with a codimension $c$ submanifold, and $c>2$, then replacement is always possible in the topological and PL locally flat settings\footnote{For a discussion of the codimension two situation, see \cite{CS2,CS3}. As the homotopy type of the top stratum can change in these results, and, for the non-locally flat results, even the number of strata can change, these results do not fit into the framework we are discussing here. }. The key technical reason for such replacement is the stability theorem (and its topological analogue) that says that the map of classifying spaces, 
\[
G_c/PL_c \to G/PL,
\]
is a homotopy equivalence for $c>2$. Moreover, the replacement of the submanifold $M$ can be achieved without altering the manifold $W$. 
 
If we move on to the setting of group actions (even just on manifolds), and we consider the issue of replacing fixed sets, then the situation is more subtle and not directly related to stability properties of classifying spaces. For instance, when the group $S^1$ is acting semi-freely on a manifold $X$, and $M$ is a codimension $4$ component of the fixed set of the action, then the (equivariant) replacement problem boils down to one for the submanifold $M$ in $X/S^1$. Since the neighborhood of $M$ has quotient a manifold, the theorem of Browder et al, yields replacement. In \cite{CW1} this was analyzed and it was shown that replacing the fixed set can force a change in the structure of the global manifold acted upon. 

Moreover, in \cite{CW2} it is proven that replacement of a component of the fixed set of a semi-free circle action is generally possible iff the codimension of that component is $0$ mod $4$. When the codimension is $2$ mod $4$, then the restriction map ${\mc S}^{S^1}(M\times D(\rho), \text{rel }\pa )\to {\mc S}(M)$ is in fact trivial, and the fixed set is thus rigidly determined by the action on the complement.   

It is easy to see that if $F'$ is $h$-cobordant to $F$, then there is a $G$-action on a manifold homotopy equivalent to $M$ with $F'$ as fixed set. However, this construction does not produce a manifold equivariantly simple homotopy equivalent to the original action, and is therefore not a replacement in the sense of this paper. 
 
We use the standard notation that for an orthogonal or unitary representation $\rho$, its unit disk is denoted by $D(\rho)$, and its unit sphere by $S(\rho)$ (or $D(V)$, $S(V)$ for the underlying vector space $V$). Also following standard usage, the structure set of a stratified space will be denoted by ${\mc S}(X)$. Indeed, this symbol will often denote the spectrum whose $0$-th homotopy group is that set. (See \cite{We1} for definitions and the spectrum structure.) The symbol ${\mc S}^G(M)$ is the same as ${\mc S}(M/G)$ and has an isovariant interpretation. We expect that this overuse of the letter $S$ (and ${\mc S}$) will cause no more confusion in this paper than it does throughout the literature.  

\medskip
 
\noindent{\bf Note}: Throughout this paper we will make the blanket assumption that all group actions are orientation preserving and all fixed sets of all subgroups are orientable. We also assume that all groups described as acting are nontrivial. 

\medskip  

Our main new result is that this replacement property for fixed sets of semifree circle actions when the normal representation has dimension $0$ mod $4$ is true in considerable generality for arbitrary compact Lie groups (including, of course, all finite groups).

\bigskip

\noindent {\bf Main Theorem} 
{\it Suppose that $G$ is a compact Lie group acting locally linearly on a 
topological manifold $W$. Suppose that near a $1$-skeleton of the fixed set $F$, the $G$-action can be identified with a complex $G$-bundle whose normal representation is a multiple of $2$. Then the forgetful map ${\mc S}^G(W)\to {\mc S}(F)$ is a split surjective map. In particular, for such actions, it is always possible to replace the fixed set by any simple homotopy equivalent (homology) manifold.}

\bigskip

\noindent {\bf Remark 1}: If the replaced fixed set is a manifold, we can arrange for the new action to be locally linear and with the same normal representation.

\medskip 
 
\noindent {\bf Remark 2}: This theorem also holds in the PL locally linear category. 

\medskip
 
The meaning of the $1$-skeleton condition is as follows. If an action is locally linear, then by definition, there is a representation $\rho$ normal to the fixed set; in other words, the equivariant tangential germ data is linearized at a $0$-skeleton. (It is only well-defined up to (stable) topological equivalence, of course.) The hypothesis we make is somewhat stronger than that $\rho=2\rho'$ for some complex representation $\rho'$. Essentially the issue is something like orientability, a condition on the germ about a $1$-skeleton: When the identification of the normal structure with a representation is transported as we move around the fixed set, do we obtain ``nontrivial monodromy''? Our complexity hypothesis ensures that we do not. In Theorem 2, stated in section 2, this condition on monodromy will be weakened. 

In the smooth case, the structure group of the equivariant normal bundle reduces to a product of compact groups of type $O$, $U$, and $Sp$ according to whether the irreducible summands of $\rho'$ are of type ${\bb R}$, ${\bb C}$, and ${\bb H}$. Only type ${\bb R}$ would allow nontrivial monodromy. However, even for smooth $G$-manifolds, our results will only produce topological $G$-manifolds.   
 
In the special case of odd order abelian groups, this theorem was proved in \cite{CW2} making use, in part, of the constructions from \cite{WY1}. At the time we had speculated that something like the above might be true. In the intervening time we verified many special cases by complicated ad-hoc constructions. The general and simple proof\footnote{However, the brevity of this paper is facilitated by references to the theory developed in \cite{CW2}.} given in section 4 below came as a pleasant surprise to us. Essentially it reduces the general case to that of unitary groups, but not by means of any fixed embedding of the group in question into the unitary group (as occurs, for instance, in the ``holomorphic induction'' proof of Bott periodicity). In this way, our main theorem is, in essence, deduced from the theorem of Browder et al referred to above.   
 
In section 2, we will review some results of \cite{CW2}, and will give some new examples, in particular showing strong differences between the PL locally linear and topological locally linear categories. In section 3, we review the methods of \cite{CW2} and the technique of stratified $\pi$-$\pi$ structure. The proofs will be given in section 4. 
 
Section 5 will give an interesting example of stratified product, namely that for $n>1$, there is a stratified space ${\bb C}^{2n}\cup D^3$ that induces an isomorphic map $L_k(\pi) \to L_{k+4n}(\pi)$. Moreover, there is a straightforward extension where $\pi$ is replaced by a general stratified space $X$. The space ${\bb C}^{2n}\cup D^3$ arises in the proofs of our main theorem. The appearance of $D^3$ essentially accounts for the role of the embedding theorem of Browder et al, and that this is an isomorphism explains why replacement of $2\times$complex representations can be reduced to (a proof of) that classical 
theorem. 
  
For $n=1$, this looks similar to the classical periodicity, but actually we do not have such a space for $n=1$ (per se; there is an appropriate example of a $3$-cell complex in \cite{WY1}, but it has rather different features). For $n>1$, the stratified space looks far from having a standard signature equal to $1$ as in the classical periodicity of Wall and Siebenmann, in that its middle dimensional homology vanishes (indeed, it is a homotopy sphere), so that the result is perhaps somewhat unexpected.   
 
We are happy to contribute this paper to the special issue of the Pure and Applied Quarterly in honor of Tom Farrell and Lowell Jones. As we shall see below, two ideas from their seminal paper \cite{FJ} play a role in our story.

\section{Review of Previous Results}

As we mentioned in the introduction, the most classical approach to replacement theorems is by analyzing the classifying spaces of neighborhoods. The earlier paper \cite{CW2} developed a rather different perspective that showed that these problems can be intimately tied to issues regarding product formulae in stratified surgery. This perspective is the one we adopt here, although our approach to the relevant product formulae is entirely different and non-computational. For the convenience of the reader, 
in this and the next sections we summarize these ideas and their earlier application. We also take advantage of this opportunity to correct some computational errors made in some examples, in the addendum to Theorem 0.1 and to Theorem 2.5, of that paper. 
 
The main previous positive results regarding replacement of fixed point sets are the following \cite{CW2}.

\bigskip

\noindent {\bf Theorem 1}
{\it If $G$ is a finite group of odd order acting on a manifold, such that the small gap hypothesis holds (i.e., no stratum is codimension two in another), and the action is smoothable near a $1$-skeleton, then one can replace the fixed point sets. Moreover, the new action can be assumed to be on the original manifold.} 

\bigskip
 
Note, though, that there is no requirement that the normal representation be twice a complex one. It is automatically complex, since the group is odd order, but replacement is true without the doubling. In contrast, the doubling is necessary for the group $S^1$, as mentioned in the introduction. 

The next result shows that the condition on monodromy (i.e. the condition regarding a neighborhood of a 1-skeleton) can be dropped in some special cases.

\bigskip

\noindent {\bf Theorem 2}
{\it If $G$ is cyclic, the action satisfies the small gap hypothesis, and the normal representation is $2\rho'$ for some complex representation $\rho'$, then replacement is possible. For general $G$, the same holds assuming in addition that the representation satisfies the strong gap hypothesis. }

\bigskip

The strong gap hypothesis (i.e., that any stratum is less than half the dimension of any larger stratum that contains it) implies (at the stratified homotopy level) the triviality of the monodromy by a theorem of Browder \cite{Sc}\footnote{Recall our blanket assumption about the orientability of all fixed sets.}. Then the replacement follows from the main theorem.

For cyclic group actions, a stabilization/destabilization procedure can be used. The idea is to view the replacement of the fixed set of a $G$-action is the same as splitting the restriction map ${\mc S}^G(M)\to {\mc S}(F)$. Moreover, the strong replacement means the splitting is trivial when composed with the forgetful map ${\mc S}^G(M)\to {\mc S}(M)$. Consider the diagram
\[
\begin{CD}
{\mc S}(M) @<<<  {\mc S}^G(M)  @>>>  {\mc S}(F) \\
@V{\cong}VV  @V{\cong}VV  @| \\ 
{\mc S}(M\times D(V), \text{rel }\pa) @<<< {\mc S}^G(M\times D(V), \text{rel }\pa) @>>> {\mc S}(F) 
\end{CD}
\]
In the diagram, we choose a suitable representation $V$ such that $M\times V$ satisfies the strong gap hypothesis, and that $V$ induces isomorphisms on the two vertical arrows. Such representations are ``periodicity representations'' and exist for cyclic groups (also given by twice of some complex representation). They were first constructed in \cite{Y, Y2}, and they exist in the generality needed here by \cite{WY1, WY2}. The commutativity of the diagram is part of the general theory.

If one has replacement under the strong gap hypothesis, then on the bottom line one has a splitting. By the periodicity isomorphism, the top line splits as well. This, combined with the main theorem, proves Theorem 2 for the cyclic group case.
 
Theorem 2 corrects the statement in the addendum to Theorem 0.1 in \cite{CW2} asserting the above result for all abelian groups. Presumably that is true, but our current methods only yield the full result for cyclic groups. We had miscalculated that in the representation theory of odd order groups, it would always be possible to find the suitable ``periodicity representation'' $V$ to achieve the strong gap hypothesis. Unfortunately, this is not even true for most representations of $G = {\bb Z}_p\times {\bb Z}_p$. 
 
We note another difference between Theorem 1 and Theorem 2. In the latter, it is not possible to assume that the new action is on the original manifold. There are rigidity theorems.

For the case of semi-free actions with simply connected fixed set, we have more or less complete answer to the replacement problem. Recall that a compact Lie group acts semi-freely and locally linearly on a manifold if and only if it acts freely and linearly on a sphere. Such group must be either $S^1$, or $S^3$, or finite with the 2-sylow subgroup being either cyclic or quaternionic. Moreover, the action by $S^1$ on the normal fibre is the complex multiplication on ${\bb C}^n$, and the action by $S^3$ is the quaternion multiplication on ${\bb H}^n$. If the group is finite with cyclic 2-sylow subgroup, then the dimension of the normal fibre is even. If the group is finite with quaternionic 2-sylow subgroup, then the dimension of the normal fibre is a multiple of $4$. 

\bigskip

\noindent {\bf Theorem 3}
{\it Suppose that $G$ is a compact Lie group (not necessarily connected) acting semi-freely and PL locally linearly on a manifold with simply connected fixed set $F$ of codimension other than $2$. Then it is possible to replace $F$ by $F'\in {\mc S}(F)$ in such a way that the complement of $F$ is unchanged, if and only if we are in one of the following cases:
\begin{enumerate}
\item $G$ is connected, and either the normal representation is $2\rho'$ for some complex representation $\rho'$, or $F'=F$;
\item $G$ is finite with cyclic $2$-sylow subgroup, and either the codimension is $0$ mod $4$, or the codimension is $2$ mod $4$ and the Kervaire classes of $F'$ in $H^{4i+2}(F,{\bb Z}_2)$ all vanish;
\item $G$ is finite with quaternionic $2$-sylow subgroup, and either the codimension is $0$ mod $8$, or the codimension is $4$ mod $8$ and the Kervaire classes of $F'$ in $H^{4i+2}(F,{\bb Z}_2)$ all vanish. 
\end{enumerate}
}

\bigskip

Recall that normal invariants are classified by homotopy classes of maps $[F,G/Top]$, and have Kervaire characteristic classes $k^{4i+2}\in H^{4i+2}(F,{\bb Z}_2)$ (see \cite{RS}). The theorem illustrates a very partial rigidity for even order groups that does not happen for odd order groups. However, such rigidity does not go as far as the rigidity that is present for the positive dimensional case.  

Theorems 3 corrects the quaternionic statement in Theorem 2.5 of \cite{CW2}. It is interesting to compare this with the analogous theorem for the topological locally linear category.

\bigskip

\noindent {\bf Theorem 4}
{\it Suppose that $G$ is a compact Lie group (not necessarily connected)  acting semi-freely and topologically locally linearly on a manifold with simply connected fixed set $F$ of codimension other than $2$. Moreover,  suppose that $H^*(F,{\bb Z})$ has no $2$-torsion. Then it is possible to replace $F$ by $F'\in {\mc S}(F)$ in such a way that the complement of $F$ is unchanged, if and only if we are in one of the following cases:
\begin{enumerate}
\item $G$ is connected, and either the normal representation is $2\rho'$ for some complex representation $\rho'$, or $F'=F$;
\item $G$ is finite with cyclic $2$-sylow subgroup, and either the codimension is $0$ mod $4$, or the codimension is $2$ mod $4$ and the Kervaire classes of $F'$ in $H^{4i+2}(F,{\bb Z}_2)$ all vanish;
\item $G$ is finite with quaternionic $2$-sylow subgroup. 
\end{enumerate}
} 

We wonder whether the topological replacement is always possible for all semi-free topologically locally linear quaternionic group actions. 

Finally, we remark that the strong replacement for semi-free PL locally linear actions is completely analyzed by Theorem 0.2 of \cite{CW2}.

\section{Review of Previous Methods}

The main idea of \cite{CW2} is to view the replacement problem as one of changing the base of a block bundle. That paper was written in the PL locally linear setting, but it is straightforward to rewrite the paper in the language of stratified surgery theory \cite{We1} to cover the topological setting. Not surprisingly, the critical issues remain unchanged from the ad-hoc approach used in \cite{CW2}. But the contrasting Theorems 3 and 4 show that there is a difference in the resulting detailed calculations. 
 
For definiteness, we outline the idea in the PL language. 

The regular neighborhood $N$ of the fixed set $F=M^G$ can be viewed as a $G$-block bundle $N\to F$. Given $F'\in{\mc S}(F)$, the problem of changing the base from $F$ to $F'$ is to make the induced map $\pa N\to F\simeq F'$ on the boundary of the neighborhood homotopic to a $G$-block bundle over $F'$. Equivariantly taking the mapping cylinder of this block bundle and gluing it to the complement $\overline{M-N}$ of the interior of the regular neighborhood $N$, produces the new $G$-manifold $M'=\overline{M-N}\cup \pa N\times[0,1]\cup F$ in ${\mc S}^G(M)$ that restricts to $F'\in {\mc S}(F)$. 

We note that by taking the quotient of the group action, the problem of changing the base of the block bundle $\pa N\to F$ (which has a $G$-sphere as fibre) is the same as the problem of changing the base of the block bundle $\pa N/G\to F$ (which has a stratified space as fibre). ``Strong replacement'', which involves knowing that the new equivariant manifold is the original one (i.e., $M'=M$), requires changing the base of the ``bubble quotient'' block bundle $N\cup_{\pa N}\pa N/G\to F$ introduced in \cite{CW2}, but which will not be discussed here. 
 
What is the obstruction to the problem of changing the base from $F$ to $F'$? Notice that the structure set ${\mc S}(F)$ can be viewed as an obstruction to a very simple base change problem: Given a homotopy equivalence $F\simeq F'$, the problem of changing the base of the identity block bundle $F\to F$ to $F'$ is the same as making the map $F\to F'$ homotopy equivalent to a block bundle map. Since for a map $F\to F'$ to be a block bundle map is the same as for the map to be a homeomorphism, the obstruction to changing the base of the identity block bundle from $F$ to $F'$ may be identified with the element that $F'$ represents in ${\mc S}(F)$. Using this, changing the base of the block bundle $\pa N/G \to F$ to $F'$ is the same as the transfer of the element $F'\in {\mc S}(F)$ to the total structure space ${\mc S}(\pa N/G \to F)$ of the block bundle $\pa N/G \to F$. 

A geometric argument in \cite{CW2} shows that the transfer ${\mc S}(F)\to {\mc S}(\pa N/G \to F)$ can be computed by knowing the ``stratified homotopy type'' of the block bundle $\pa N/G \to F$ over a $1$-skeleton. (This argument uses a thickening ${\bb C}P^2\times F$ of $F$ to ensure enough of a skeleton over each simplex of $F$. It is reminiscent of the transfer trick of \cite{FJ} although used for a completely different reason.) Therefore if the block bundle is trivial over a $1$-skeleton, then the transfer is the same as the map 
\[
\times S(\rho)/G\colon {\mc S}(F)\to {\mc S}(F\times S(\rho)/G \to F)
\]
induced by product with the fibre, where $\rho$ is the normal $G$-representation of the fixed set $F$. Then as explained in \cite{CW2}, by the surgery theory that computes the structure sets, the vanishing of the product map is a consequence of the vanishing of the map 
\[
\times S(\rho)/G\colon L(?)\to L^{BQ}(?\times S(\rho)/G)
\]
induced by product on stratified $L$-groups. For the strong replacement, the problem is a consequence of the vanishing of the product with the bubble quotient
\[
\times (D(\rho)\cup S(\rho)/G)\colon L(?)\to L^{BQ}(?\times (D(\rho)\cup S(\rho)/G)).
\]

The stratified $L$-groups used here are the ones introduced by Browder and Quinn \cite{BQ} for transverse isovariant surgery. In \cite{We1}, it is shown that these $L^{BQ}$ play a key role in the non-transverse theory as well. Indeed, in that theory, they describe both the global surgery obstruction, and in a cosheaf form, the normal invariants.

The method we use to prove the vanishing of the product map on the stratified $L$-groups is rather geometrical, and makes substantial use of stratified spaces. This avoids most of the hard computation. The general setting for this, as used in \cite{WY1,WY2,Y,Y2}, is a stratified space $\Sigma$ with the property that $L^{BQ}(X\times\Sigma) = 0$ for any $X$. Then for any union $\Sigma'$ of some closed strata of $\Sigma$, the product 
\[
\times \Sigma'\colon L(X)\to L^{BQ}(X\times\Sigma')
\]
vanishes because it is a composition
\[
L(X)\xrightarrow{\times \Sigma} L^{BQ}(X\times\Sigma)=0
\xrightarrow{\text{restrict}} L^{BQ}(X\times\Sigma'). 
\]
If we can express $S(\rho)/G$ as the union of some closed strata for suitable $\Sigma$, then we proved the replacement theorem.

How do we find suitable $\Sigma$? Our idea comes from Wall's $\pi$-$\pi$ theorem (see \cite{Wa1} Theorem 3.3): If $M$ is a connected manifold with connected boundary $\pa M$, such that the map $\pi_1(\pa M)\to \pi_1(M)$ induced by the inclusion is an isomorphism, then the surgery obstruction group $L(M)$ of the pair $(M,\pa M)$ vanishes. The theorem can be extended to stratified spaces with $\pi$-$\pi$ structures. The equivariant (actually isovariant) version of such stratified $\pi$-$\pi$ structures is introduced in Definition 1.1 of \cite{Y2}. By Propositions 1.2 and 1.4 of \cite{Y2}, if $\Sigma$ has stratified $\pi$-$\pi$ structure, then $L^{BQ}(X\times\Sigma)=0$ for any $X$.

Suppose $\Sigma=\hat{\Sigma}/G$ for a homotopically $G$-stratified space $\hat{\Sigma}$. Then it is easy to see that $\Sigma$ has stratifed $\pi$-$\pi$ structure if and only if $\hat{\Sigma}$ has $G$-isovariant $\pi$-$\pi$ structure. Now for a union $\hat{\Sigma}'$ of some closed strata of $\hat{\Sigma}$ and any subgroup $H\sub G$, we also know that the product
\[
\times \hat{\Sigma}'/H\colon L(X)\to L^{BQ}(X\times\hat{\Sigma}'/H)
\]
vanishes because it is a composition
\[
L(X)\xrightarrow{\times \hat{\Sigma}'/G=0} L^{BQ}(X\times \hat{\Sigma}'/G)
\xrightarrow{\text{forget}} L^{BQ}(X\times\hat{\Sigma}'/H),
\]
where the forgetful map simply considers $G$-spaces as $H$-spaces. This shows that the vanishing of the product on the stratified surgery obstruction can descend from a group action to a subgroup action. The fact will be used in the proof of the main theorem, and is consistent with the following fact (a tautology by the definition of replacement): Suppose that $H\sub G$ is a subgroup, and $\rho$ is a unitary representation of $G$ such that only the origin is fixed by $H$. If fixed points with $\rho$ as normal representation has (strong) replacement for $G$, then $\rho$ also has (strong) replacement as an $H$-representation.

Next we explain how the general ideas above can be applied to specific cases. 

When we are dealing with finite group actions (or locally free in the sense that all isotropy groups are finite), away from the prime 2 all of the $L^{BQ}$ that occur in stratified surgery naturally split. Consequently, the isovariant structure set splits into pieces that are concentrated on the components of the singularity set. This is because using ${\bb Q}$ as the coefficient ring in place of ${\bb Z}$, permutation modules are all projective, so that one can use the singular chain complexes to give a description of the relevant $L$-groups. On the other hand, for any $X$, an induction beginning with Ranicki's localization theorem \cite{R} for Wall's surgery obstruction groups shows that $L^{BQ}(X) \to L^{BQ}(X,{\bb Q})$ is an equivalence after $\otimes{\bb Z}[1/2]$.

At the prime $2$, for a representation $\rho$ of an odd order group $G$, L\"{u}ck and Madsen \cite{LM} showed that some odd multiple of $S(\rho)/G$ is the boundary of a manifold in essentially the $\pi$-$\pi$ fashion. By taking $\Sigma$ to be the manifold with (certain odd multiple of) $S(\rho)/G$ as the boundary and taking $\Sigma'$ to be $S(\rho)/G$, this implies that for odd order groups, the surgery obstruction splits even integrally. Such integral splitting lies behind the proof of Theorem 1 in \cite{CW2}.

The main theorem is proved by identifying $S(\rho)/G$ as the union of some closed strata in a stratified space with $\pi$-$\pi$ structure. Unlike the odd order group actions, here we do not expect to find suitable manifolds with (perhaps odd multiples of) $S(\rho)$ as boundary, so that the use of stratified spaces is essential. (See \cite{DS} Proposition 3.7 for a discussion of the difficulties in finding such $G$-manifolds for $G = {\bb Z}_2$.) For abelian groups, such stratified spaces are constructed in \cite{WY1} explicitly for the case $\rho=2\rho'$ with $\rho'$ being any irreducible representation. Then the special cases may be combined to obtain the suitable stratified space for the general $\rho=2\rho'$, with growing complexity in the stratification: There were $3\dim_{\bb C}\rho'-1$ strata.

For the group $S^1$ acting on ${\bb C}^k$ by scalar multiplication, the replacement is computed by the effect of the product with $S({\bb C}^k)/S^1={\bb C}P^{k-1}$ on surgery obstructions. It is then clear why one has a dichotomy between replacement and rigidity depending on whether one has a ${\bb C}P^{\rm odd}$, which kills surgery obstructions, or a ${\bb C}P^{\rm even}$, which gives rise to periodicity isomorphisms on surgery groups. For the group $S^3$ acting on ${\bb H}^k$ by quaternion multiplication, we need to consider the product with quaternionic projective spaces, and have similar dichotomy.  (For general spheres of representations, of course, one tends to neither kill all nor preserve all surgery obstructions, leading to a more varied range of phenomena in the replacement problem.) 
 
The original purpose of introducing stratified spaces with $\pi$-$\pi$ structure was to contruct isovariant periodicity maps ${\mc S}^G(M) \to {\mc S}^G(M\times D(\rho), \text{rel }\pa)$ for suitable $G$-representations $\rho$. Unfortunately, the ad-hoc techniques used for particular groups and particular representations did not suffice for general Lie groups. In \cite{WY2} the authors constructed a much more efficient stratified space for representation spheres (or better, an improved construction of periodicity spaces) using a key calculation at the core of \cite{FJ} -- the symmetric square construction.    
 
What the argument in \cite{WY2} directly shows, however, is that a representation of the form $\rho\oplus\epsilon^4$ is a periodicity representation. For the periodicity problem, one can, after the fact remove the trivial summand $\epsilon^4$. But for the replacement problem this cannot be done: One would only obtain the result that $F'$ can be isovariantly embedded in $F\times D^4$, which was obvious anyway as $F'$ embeds in $F\times D^4$. While we seemed no closer to proving that the general $\rho=2\rho'$ is a ``replacement representation'', it then became a more pressing issue to try to prove such a replacement theorem. This will be accomplished in the following section.

\section{Proofs of Theorems}

We begin with a geometrical observation (see \cite{Da} for a much more complete analysis).

\bigskip
 
\noindent {\bf Lemma} 
{\it Let $U(n)$ act on ${\bb C}^{2n}$ as twice the defining action on ${\bb C}^n$. Then $S({\bb C}^{2n})/U(n) \cong D^3$ for $n>1$.}

\bigskip

\begin{proof}
The action on the sphere has two orbit types, depending on whether the pair of vectors $(v,w)\in S({\bb C}^{2n})$ span a $1$- or $2$-dimensional subspace of ${\bb C}^n$. If they span a $1$-dimensional subspace, then the ratio between $v$ and $w$ gives a homogeneous coordinate for a point in the complex projective line ${\bb C}P^1\cong \pa D^3$. Note that although the complex line spanned by the vectors changes as one moves around the $U(n)$-orbit of the pair, the projective element represented by the pair within that complex line is well defined.
 
In general, we decompose $w$ into projections in the direction of $v$ and in the direction orthogonal to $v$:
\[
w=\proj_vw+\proj_{v^{\perp}}w.
\]
The $U(n)$-orbit of the pair uniquely corresponds to the pair $([v,\proj_vw],\|\proj_{v^{\perp}}w\|)$. The first coordinate $[v,\proj_vw]\in {\bb C}P^1\cong \pa D^3$ is the point on the complex projective line given by the ratio between $v$ and $\proj_vw$. The second coordinate $\|\proj_{v^{\perp}}w\|$ describes the radial coordinate of the ball $D^3$. As a matter of fact, it is more natural to take $1-\|\proj_{v^{\perp}}w\|$ as the radius. The case $\|\proj_{v^{\perp}}w\|=1$ corresponds to the pair $(0,w)$. Since all such pairs lie in the same $U(n)$-orbit, the case contains only one point, and gives the origin of the ball. On the other hand, the case $\|\proj_{v^{\perp}}w\|=0$ means that $v$ and $w$ span a $1$-dimensional subspace of ${\bb C}^n$. This corresponds to the boundary of the disk $D^3$ and is taken as the lower stratum of the stratified space $D^3$.

\end{proof}

\begin{proof}[Proof of Main Theorem]
As explained in section 3, the proof of the main theorem may be reduced to showing that the product with $S(\rho)/G$ induces the trivial map on the stratified $L$-groups. For the special case that $G=U(n)$ and $V={\bb C}^{2n}$ as in the lemma, the product 
\[
\times S({\bb C}^{2n})/U(n) =\times D^3\colon
L(F)\to 
L^{BQ}(F\times D^3)
\]
is trivial because the stratified $L$-group $L^{BQ}(F\times D^3)$ is trivial. Here the triviality of $L^{BQ}(F\times D^3)$ is obtained by applying the $\pi$-$\pi$ theorem to the pair $F\times (D^3,S^2)$.

Under the assumption of the main theorem, we need to show that the product $\times S(2\rho')/G$ induces the trivial map on stratified $L$-groups. The $G$-representation $\rho'$ may be embedded into the defining representation of $U(n)$ on ${\bb C}^n$ (via a homomorphism $G\to U(n)$). As explained in section 3, the triviality of the map $\times S({\bb C}^{2n})/U(n)$ on the stratified $L$-groups implies (or descends to) the triviality of the map $\times S(2\rho')/G$ on the corresponding stratified $L$-groups.

The main theorem is then proved if the dimensions of all the strata in the quotient that touch the fixed set are of dimension at least five, so that we can apply surgery theory. Note, though, that the dimension of the fixed set is at least $3$ (for there to be any other manifold homotopy equivalent to the fixed set). The condition on $\rho$ guarantees that the quotient of the sphere by the action is at least $2$-dimensional (with the extreme case being $2\times$the defining representation for the circle), so there are no low dimensional complications. 

\end{proof}

The proof above of the main theorem is based on the $\pi$-$\pi$ structure given in lemma. The argument gives a conceptual a-priori explanation for the replacement part of the Browder-Casson-Haefliger-Sullivan-Wall embedding theorem. (See \cite{Wa1} Corollary 11.3.1.) The book \cite{We1} contains three proofs of the full result. 

Now we turn to the proofs of the two remarks to the main theorem. 
 
Remark 1 asserts that when the replaced fixed set is assumed to be a manifold, then one has local linearity of the replacement. This follows from a relative form of the main theorem.   

The manifolds $F'$ simple homotopy equivalent to $F$ are the same as those simple homotopy equivalent to $\mathring{F}$ ($F$ punctured) rel $\pa$ by gluing in the ``missing'' disk. Applying the relative version of the main theorem to $\mathring{M}$, and gluing back in the disk of $M$ that we removed, gives an equivariant manifold $M'$ with fixed set $F'$ that is locally linear near the center of this disk. This implies, in the topological category, local linearity by the topological homogeneity of homotopically stratified spaces (see \cite{Q} Corollary 1.3). 
 
Remark 2 regarding the PL locally linear category requires no explanation except for one low dimensional issue. Whereas the topological category uses stratified and controlled surgery, the PL case uses instead blocked surgery, whose formal obstruction theory is identical aside from the decorations for the $L$-groups. However, blocked surgery requires that all blocks have no low dimensional strata, which would naively restrict the (dimensions of strata in the) normal representations to which our theorem could apply. 
 
Using Theorem 12.1 of \cite{Wa1}, we can arrange for any simple homotopy equivalence $F'\to F$ to be a PL homeomorphism near a $1$-skeleton (see proof of Theorem 1.7 in \cite{CW2}). Then the blocked surgery starts over simplices of dimension $2$. Since the strata in stratified spaces $S(\rho)/G$ used in our proof have dimension at least $3$, the corresponding blocks over simplices of dimension $\ge 2$ have dimension at least $5$. Therefore there is really no restriction on the strata of the normal representation,  completing the argument for Remark 2.

Now we turn to the proof of Theorems 3 and 4. 
 
\begin{proof}[Proof of Theorem 3]
Similar to the proof of the main theorem, the key issue is the computation of the transfer $\times S(\rho)/G$ on the surgery obstructions. In case $G=S^1$, we have $\rho={\bb C}^k$, $S(\rho)/G={\bb C}P^{k-1}$, so that the product $\times S(\rho)/G$ is trivial on the surgery obstruction when $k$ is even and is isomorphic on the surgery obstruction when $k$ is odd. This implies the replacement in case $k$ is even (so that $\rho=2\rho'$ for $\rho'={\bb C}^{k/2}$) and rigidity in case $k$ is odd. In case $G=S^3$, we have $\rho={\bb H}^k$, $S(\rho)/G={\bb H}P^{k-1}$, and the argument is similar.

The case of finite $G$ is more complicated because the products are sometimes neither trivial (replacement) nor an isomorphism (rigidity). Now $S(\rho)/G$ is a space form, which represents the trivial element in the symmetric $L$-group $L_n({\bb Z}G)\otimes{\bb Z}[1/2]$ (as mentioned earlier in our discussion of why $L^{BQ}$ splits). Thus we can work localized at the prime $2$, where the normal invariants are determined by characteristic classes. Moreover, we can therefore, following \cite{Wa2}, restrict our attention to the cover $S(\rho)/G(2)$ corresponding to the $2$-sylow subgroup $G(2)$, which are either lens spaces or quaternionic space forms.

According to \cite{HMTW}, the product $\times S(\rho)/G(2)$ on the surgery obstruction is determined by applying ``$\kappa$-homomorphisms'' to mod $2$ reductions of characteristic classes, which are a combination of Wu class and the Morgan-Sullivan $L$-class \cite{MS}. Reduced mod $2$, the latter is the square of the Wu class, so that everything is determined by Wu class. By the Wu formula for Wu class, these are determined via the action of the Steenrod algebra on mod $2$ cohomology, in any case a $2$-local homotopy invariant. Note that the product on the surgery obstruction depends only on the homotopy class of the space form $S(\rho)/G(2)$, and the homotopy type of space forms of a given dimension is classified by their $k$-invariants. Since the $k$-invariants of space forms must be odd (by, for instance, the Borsuk-Ulam theorem), the product with space forms depends only on the dimension, as far as the $2$-localized surgery obstructions are concerned.

Since only the dimension of $\rho$ really matters, we only need to compute the transfer for one conveniently chosen space form. For the case $G(2)$ is cyclic, we may choose $\rho=k\tau$, where $\tau={\bb C}$ with the cyclic group action given by complex multiplication. For the case $G(2)$ is quaternionic, we may choose $\rho=k\tau$, where $\tau={\bb H}$ with the quaternion group action given by quaternion multiplication. If $k$ is even, then this gives the vanishing of the transfer in dimensions a multiple of $4$ in the cyclic case and a multiples of $8$ in the quaternionic case. If $k$ is odd, then the product of $\times S(\rho)/G(2)$ with the Kervaire surgery problem (generating $L_2(e)$) is nontrivial, and these (and their obvious descendents) are the only nontrivial products on $L$-theory with coefficients (see \cite{MS}).

These tell us the homological part of the obstruction to base change in the block fibrations\footnote{Note that the obstruction to base change is in the fiber of an assembly map. By the homological part of the obstruction we mean the part that lies in the homology, and when this vanishes, one has another second order obstruction in light of the exact sequence of a fibration. }. There is, in general, a final surgery obstruction. It is at this point that we use the hypothesis of simple connectivity, and a variant of the trick used in the proof of Remark 1 to complete the proof. (Since we are working in PL rather than Top, we cannot appeal to \cite{Q}.) We puncture the fixed set, and replace that incurring no final normal cobordism. At that point we cone the boundary, and obtain a new group action that is PL and locally linear except perhaps at the final cone point. However, the obstruction to local linearity there is determined by an element of (the Tate cohomology of) the Whitehead group, see \cite{We3} (or \cite{CW3} if working up to concordance, which is adequate for us). Since we are working throughout with simple homotopy equivalences, the vanishing is a tautology. 

\end{proof}

\begin{proof}[Proof of Theorem 4]
For the topologically locally linear actions, since the stratified space in the lemma actually has two strata, the teardrop neighborhood theorem of \cite{HTWW} directly applies, and the vanishing of the relevant surgery obstructions happens in $L^{-f}$ (the negative $L$-group). A fortiori, calculation in $L^p$ suffices. Vanishing here on $L_2(e)$ is well known (see \cite{HMTW}) but to complete the map of spectra $\times S(\rho)/G(2)\colon L(e)\to L(G(2))$ would require knowing vanishing for homotopy with coefficients, and in particular computing the product map $L_3(e,{\bb Z}_2) \to L_6(G(2),{\bb Z}_2)$.  However, if the fixed set $F$ does not have torsion in its cohomology, the information about integral homotopy is enough and the proof is completed in the same way as Theorem 3. 

\end{proof}

The replacement problem for non-semi-free actions are much more complicated. In the remaining part of the section, we carry out some preliminary computations related to the actions of $S^3=SU(2)$.

Recall that the irreducible representations of $SU(2)$ can be described in terms of spaces of polynomials in two variables. Denote by $\rho_d$ the irreducible complex representation of $SU(2)$ acting on the space of degree $d$ polynomials on ${\bb C}^2={\bb C}x\oplus {\bb C}y$. The complex dimension of $\rho_d$ is $d+1$. We will study the replacement problem when the normal representation of the fixed set is $\rho=k\rho_d$.
 
The maximal torus $T$ of $SU(2)$ is a circle and is (up to conjugation) the subgroup consisting of the actions
\[
x\to \lambda x,\s y\mapsto \lambda^{-1} y,\s |\lambda|=1.
\]
The eigenspaces for $T$ are given by monomials, and the Weyl group is ${\bb Z}_2$ with a representative of the nontrivial element given by the action (of order $4$ in $SU(2)$)
\[
x\mapsto -y,\s y\mapsto x.
\]
 
The representation $\rho_d$ behaves quite differently for even and odd $d$. For even $d$, the center $-I$ acts trivially (so the action is not effective), and the monomial $(xy)^{d/2}$ is fixed by the whole $T$. For odd $d$, the action is effective, and all isotropy groups are finite (i.e., the action on the unit sphere is locally free). 
 
If $k$ is even, then the main theorem implies replacement for $\rho=k\rho_d$. If $k$ is odd, the stabilization/destabilization procedure can reduce the case $\rho=k\rho_d$, $k>3$, to the case $\rho=3\rho_d$. However, the procedure cannot be used to further reduce the case to $\rho=\rho_d$, because $k\rho_d$ satisfies the small gap hypothesis (a condition for the stabilization/destabilization procedure) if and only if $k\ge 3$. 

Since the action is not semi-free, the stratified space $S(\rho)/G$ can be very complicated. However, for some subgroup $H\sub G$, the image of the $H$-fixed set $S(\rho)^H/(NH/H)$ in the quotient $S(\rho)/G$ may be simpler (for example, in case $NH/H$ acts freely on $S(\rho)^H$). Note that the product map $\times S(\rho)^H/(NH/H)$ is a composition
\[
\times S(\rho)^H/(NH/H)\colon
L(?)\xrightarrow{\times S(\rho)/G} L(?\times S(\rho)/G)
\xrightarrow{\text{restrict}} L(?\times S(\rho)^H/(NH/H)),
\]
so that if $\times S(\rho)^H/(NH/H)$ is nontrivial on the surgery obstructions, then $\times S(\rho)/G$ is also nontrivial. Therefore we may get certain rigidity to the replacement problem by ``picking out'' the fixed sets of some subgroups.

For the case $\rho=\rho_d$, $d$ even, we take $H=T$ and get $S(\rho_d)^T=S^1$, $S(\rho_d)^T/(NT/T)=S^1/{\bb Z}_2=S^1$. Since the product with a circle is injective on structure sets (see \cite{Sh,Wa1}) aside from decorations, we see that there is rigidity (up to decorations) in this case. 

Note that for the case $\rho=k\rho_d$, $d$ even and $k>1$ odd, we have $S(\rho_d)^T/(NT/T)=S^{2k-1}/{\bb Z}_2={\bb R}P^{2k-1}$. However, the product with ${\bb R}P^{2k-1}$ is no longer injective on structure sets. This specifically illustrates why the case $\rho=k\rho_d$ for odd $k>1$ cannot be reduced to $\rho=\rho_d$.
 
Next we consider $\rho=\rho_d$, with $d=2r+1$ odd. For the cyclic subgroup $H={\bb Z}_d\sub T\sub SU(2)$, the fixed set $S(\rho_d)^{{\bb Z}_d}$ is the unit sphere of the complex vector space spanned by $x^d$ and $y^d$, and we have 
\[
\dfrac{S(\rho_d)^{{\bb Z}_d}}{N{\bb Z}_d/{\bb Z}_d}
=\dfrac{S({\bb C}x^d\oplus {\bb C}y^d)}{ax^d+by^d\sim a\lambda^dx^d+b\lambda^{-d}y^d\sim -ay^d+bx^d}
={\bb R}P^2.
\]
Since the product with ${\bb R}P^2$ preserves Kervaire invariants, it is at least necessary for the Kervaire cohomology classes to vanish for the replacement to be possible.

On the other hand, we can fairly completely analyze the obstruction away from the prime $2$. Recall that since $d$ is odd, the isotropy subgroups of $S(\rho_d)$ must be finite. Moreover, the unique element $-I$ of order $2$ in $SU(2)$ acts on $S(\rho_d)$ as the antipodal map, so that no even order subgroups can have fixed points. Hence all isotropy groups must be odd order cyclic, which lie in $T$ up to conjugacy. Consequently, the strata of $S(\rho_d)/SU(2)$ are in a one-to-one correspondence with odd integers $p$ up to $d$. Specifically, the closed stratum corresponding to $p$ is the quotient of the sphere of 
\[
V_p=\oplus_{p\text{ divides } d-2i}({\bb C}x^iy^{d-i}\oplus {\bb C}x^{d-i}y^d)
\]
by the group $N{\bb Z}_p/{\bb Z}_p$. The dimension gap between the strata are at least $4$. For $p=1$, the stratum is the top stratum $S(\rho_d)/SU(2)$, a simply connected rational homology manifold. For $p\ne 1$, the stratum is a rational homology manifold with the fundamental group ${\bb Z}_2$ acting orientation reversingly. 

As explained in section 3, away from the prime $2$, the stratified surgery obstruction $L^{BQ}(S(\rho_d)/SU(2))$ naturally splits according to the strata of $S(\rho_d)/SU(2)$. Since the $L$-group of ${\bb Z}_2$ with nontrivial orientation is $2$-torsion,  the surgery obstruction pieces corresponding to the non-top strata are trivial away from $2$. Therefore only the surgery obstruction corresponding to the top stratum remains, and we have
\[
L^{BQ}(S(\rho_d)/SU(2))\otimes{\bb Z}[1/2]
=L(e)\otimes{\bb Z}[1/2].
\]

Now we understand that for simply connected $X$, the product 
\[
\times S(\rho_d)/SU(2)\colon L(X)\mapsto L^{BQ}(X\times S(\rho_d)/SU(2))
\]
with the stratified space $S(\rho_d)/SU(2)$ is, after localizing away from $2$, the same as the product
\[
\times S(\rho_d)/SU(2)\colon L_n(e)\otimes{\bb Z}[1/2]\to L_{n+2d-2}(e)\otimes{\bb Z}[1/2]
\]
with the rational homology manifold $S(\rho_d)/SU(2)$ of dimension $2d-2=4r$. Such product depends only on the signature of the rational homology manifold $S(\rho_d)/SU(2)$, which is a rational quaternionic projective space. Since the projective space has signature $0$ or $1$ depending on whether $r$ is odd or even, one has (away from $2$) rigidity for $d = 1$ mod $4$ and no obstruction at all for $d = 3$ mod $4$. 
 
To summarize, suppose $\rho=\rho_d$ is an irreducible representation of $SU(2)$. Then for $d$ even, there is rigidity (up to decoration).  For $d$ odd, one always has Kervaire cohomology obstructions at the prime $2$. Away from $2$ and in case the fixed set is simply connected, there is no obstruction for $d = 3$ mod $4$, and there is rigidity for $d = 1$ mod $4$. It is reasonable to conjecture that for $d = 1$ mod $4$, one actually has rigidity integrally. But given how our obstructions at $2$ and away from $2$ come from such different sources, it seems that proving this could be difficult.

\section{A Product that Induces Periodicity in Surgery}

Finally, we return to discuss the periodicity isomorphism mentioned at the end of the introduction. A key idea underlying our proof of the replacement theorems is to construct a suitable stratified space $\Sigma$ with the property that $L^{BQ}(X\times \Sigma) = 0$ for any $X$. Now suppose $Z$ is a stratified space that contains such a $\Sigma$ as a union of closed strata. Then for any $X$, there is a product construction:  
\[  
L(X) \xrightarrow{\times Z}
L^{BQ}(X\times Z)  \xleftarrow[\simeq]{\text{incl}}
L^{BQ}(X\times Z,\text{rel }X\times\Sigma). 
\]
The construction is the composition of the natural map induced by taking the product with $Z$ and the inverse of the natural map induced by the inclusion. The inclusion map is an isomorphism because of the assumption $L^{BQ}(X\times \Sigma) = 0$. 

For our purposes, we will only consider the case where $\Sigma$ is the singular set of $Z$, i.e., $Z-\Sigma$ is a manifold. Then the product construction gives
\[  
L(X) \xrightarrow{\times Z}
L^{BQ}(X\times Z)  \xleftarrow[\simeq]{\text{incl}}
L(X\times \overline{Z-\Sigma},\text{rel }\pa). 
\]
Here $\overline{Z-\Sigma}$ is a closure of the top pure stratum as a manifold with boundary. (Alternatively, one can easily work with the $L$-theory with compact supports.) 

The basic example to which the product construction can be applied is $Z = {\bb C}P^2\cup D^3$ (glued along ${\bb C}P^1=S^2$) with $\Sigma=D^3$. The property $L^{BQ}(X\times \Sigma) = 0$ follows from the $\pi$-$\pi$ theorem. The resulting product contruction $L(X) \to L^{BQ}(X\times D^4,\text{rel }\pa)$ is Siebenmann's periodicity isomorphism.  
 
Now let us take $Z={\bb C}^{2n}\cup D^3$. The space is obtained by compactifying ${\bb C}^{2n}$ by attaching its unit sphere $S^{4n-1}$ at infinity and then collapsing the unit sphere to the quotient $D^3 = S^{4n-1}/U(n)$ via the lemma in section 4. This is exactly the bubble quotient construction of the biaxial unitary action. 
 
Similar to the basic example, with $\Sigma=D^3$, this $Z$ is suitable for the product construction. And our conclusion is that the construction gives an isomorphism $L(X) \to L^{BQ}(X\times D^{4n},\text{rel }S^{4n-1})$, just like Siebenmann's isomorphism.

From the present perspective this is clear: The bubble quotient $Z={\bb C}^{2n}\cup D^3=D^{4n}\cup S^{4n-1}/U(n)$ for the $U(n)$-representation ${\bb C}^{2n}$ may be compared with the bubble quotient $Z'={\bb C}^{2n}\cup {\bb C}P^{2n-1}=D^{4n}\cup S^{4n-1}/S^1$ for the restriction of the $U(n)$-representation to the subgroup $S^1$. By the discussion in section 3, the $\pi$-$\pi$ structure on $D^3=S^{4n-1}/U(n)$ descends to a $\pi$-$\pi$ structure on ${\bb C}P^{2n-1}=S^{4n-1}/S^1$, and the two product constructions can be compared. In other words, for $\rho={\bb C}^{2n}$, $G=U(n)$ and $H=S^1$, the following diagram commutes.
\[
\begin{CD}
L(X) @>{\times Z}>>  L^{BQ}(X\times (D(\rho)\cup S(\rho)/G)) @<{\text{ incl}}<\simeq<  L(X\times D(\rho),\text{rel }S(\rho)) \\
@| @VVV  @|   \\ 
L(X) @>{\times Z'}>>  L^{BQ}(X\times (D(\rho)\cup S(\rho)/H))  @<{\text{ incl}}<\simeq<  L(X\times D(\rho),\text{rel }S(\rho))
\end{CD}
\]
In particular, the product construction for the $G$-representation is isomorphic if and only if product construction for the $H$-representation is isomorphic.

If we forget stratified structure, then the bubble quotient $Z'={\bb C}P^{2n}$ is a manifold of signature $1$. Therefore the product map
\[
\times {\bb C}P^{2n}\colon 
L(X) \xrightarrow{\times Z'}  L^{BQ}(X\times Z')  
\xrightarrow{\text{forget stratification}} L(X\times {\bb C}P^{2n})
\]
is an isomorphism. On the other hand, the composition
\[
\text{incl}\colon L(X\times D^{4n},\text{rel }S^{4n-1}) \xrightarrow{\text{inclusion}}  L^{BQ}(X\times Z')  
\xrightarrow{\text{forget stratification}} L(X\times {\bb C}P^{2n})
\]
is simply the map induced by the inclusion $D^{4n}\to {\bb C}P^{2n}$ of non-stratified spaces. Since the inclusion induces an isomorphism on the fundamental groups, the composition is also an isomorphism.
 
Combining the two isomorphisms together, we see that the product construction for the $S^1$-representation ${\bb C}^{2n}$ is an isomorphism. As explained above, the product construction for the $SU(n)$-representation ${\bb C}^{2n}$ is also an isomorphism.

\end{document}